# Improved Finite Difference Method with a Compact Correction Term for Solving Poisson's Equations

Kun Zhang[1], Liangbi Wang[1] and Yuwen Zhang[2, *]

[1] College of Mechatronic Engineering, Lanzhou Jiaotong University, Lanzhou, 730070, China

[2] Department of Mechanical and Aerospace Engineering, University of Missouri, Columbia, MO 65211, USA

**ABSTRACT**

An improved finite difference method with compact correction term is proposed to solve the Poisson's equations. The compact correction term is developed by a coupled high-order compact and low-order classical finite difference formulations. The numerical solutions obtained by the classical finite difference method are considered as fundamental solutions with lower accuracy, whereas compact correction term is added into source term of classical discrete formulation to improve the accuracy of numerical solutions. The proposed method can be extended from two- to multi-dimensional cases straightforwardly. Numerical experiments are carried out to verify the accuracy and efficiency of this method.

## 1 INTRODUCTION

Poisson's equation is often encountered in mechanical engineering applications, theoretical physics and other fields. Finite-difference methods (FDM) are efficient tools for solving the partial differential equation, which works by replacing the continuous derivative operators with approximate finite differences directly [1-2]. One of the simplest and straightforward finite difference methods is the classical central finite difference method with the second-order accuracy [3]. This method is flexible to develop the discretization for solving Poisson's equation on multi-dimensional cases on uniform or non-uniform grids [4-5]. The FDM with the second order accuracy is an explicit method and it is established on a stencil of 5-points for 2D space and 7-points for 3D space. However, the significant shortcoming of the second order FDM is the solution of Poisson's equation has lower accuracy [6] and higher order accuracy method requires a larger stencil. A large stencil requires some modifications near the boundaries where the points needed in the stencil are not available, but this problems could be eliminated if a compact method is used [7-8].

Higher order compact finite difference method was first introduced by Kreiss and Oliger [9] and implemented by Hirsh [10], then popularized by Lele [11]. Compact schemes can provide numerical solutions with spectral-like resolution and very low numerical dissipation. Compared with explicit finite difference methods, this method is implicit and obtains the evaluation of derivatives with higher accuracy for the same number of grid points [12-14]. Many studies on compact finite difference schemes have been conducted for solving Poisson's equation [15-17]. Zhang [18] and Ge [19] developed a fourth-order compact scheme with multigrid method for solving Poisson's equation on two- and three-dimensional spaces, respectively. Shiferaw and Mittal [20-21] investigated the numerical solutions of three dimensional Poisson's equation with the Dirichlet's boundary conditions. The Poisson's equation was approximated by a 19-points fourth order compact finite difference approximation schemes. Lai and Tseng [22] studied a formally fourth-order accurate compact scheme for 3D Poisson's equation in cylindrical coordinates. Kyei et al. [23] derived a family of sixth-order compact finite difference schemes for the three-dimensional Poisson's equation, who considered the discretization of source function on a compact finite difference stencil. There scheme have good numerical stability and provide high accuracy approximations for solving Poisson's equations. However, the disadvantage of compact finite difference scheme is that it is an implicit form with spatial discretization and lack of flexibility due to the complex matrix transformation, especially for multi-dimensional cases on non-uniform grid.

* Corresponding author. Email: zhangyu@missouri.edu.

Therefore, it remains a challenge to develop a flexible, stable and accurate finite difference scheme for solving Poisson's equations, particularly in multi-dimensional cases and/or non-uniform grids. In fact, an improved efficient method for solving partial differential equations can be developed by combining the advantages of two different numerical methods [24-26]. Li and Li [27] studied the multigrid method combined with a fourth order compact scheme for the 2D Poisson's equation. The results showed that the new method was of higher accuracy and less computational time. Fukuchi [28] investigated finite difference method and algebraic polynomial interpolation for solving Poisson's equation over arbitrary domains. Li et al. [29-30] studied hybrid lattice Boltzmann and finite volume method for natural convection. Hejranfar and Ezzatneshan [31] proposed a new higher order compact finite difference lattice Boltzmann method for steady and unsteady incompressible flows.

The objective of this paper is to develop an improved finite difference method with compact correction term (CCFDM) for solving Poisson's equations. The proposed method has the advantage of flexibility and high accuracy by coupling high order compact and low order classical finite difference formulations. Several numerical experiments are calculated to verify the high order convergence rate of the CCFDM and the efficiency of solving Poisson's equations from two to multi-dimensional cases.

## 2 FINITE DIFFERENCE METHOD WITH COMPACT CORRECTION TERM

A two-dimensional Poisson's equation is considered first to present the basic ideas. The Poisson's equation with Dirichlet boundary condition can be written in the form of

$$\frac{\partial^2 \varphi}{\partial x^2}+\frac{\partial^2 \varphi}{\partial y^2}=f(x,y) \quad (x,y) \in \Omega \tag{1}$$

where $\Omega = [0, L_x] \times [0, L_y]$ is a rectangular domain with suitable boundary conditions defined on $\partial\Omega$. The solution $\varphi(x, y)$ and the forcing function $f(x, y)$ are assumed to be sufficiently smooth and to have the required continuous partial derivatives. The spatial domain $\Omega$ is discretized with non-uniform grid sizes $\Delta x_i = x_{i+1} - x_i$ and $\Delta y_j = y_{j+1} - y_j$ in the $x$- and $y$-coordinates, respectively. The grid points are $(x_i, y_j)$ in the $x$- and $y$-coordinate directions with the number of grid points $M \times N$, $0 \leq i \leq M$, $0 \leq j \leq N$. The quantity $\varphi(x_i, y_j)$ represents the exact solution at $(x_i, y_j)$, while $\varphi_{i,j}$ represents the numerical solution at $(x_i, y_j)$.

### 2.1 Classical and compact finite difference formulations on non-uniform grid

The classical finite difference approximation with the second order accuracy on non-uniform grid for the second order derivative in the $x$-direction may be written as follows:

$$\left.\frac{\partial^2 \varphi}{\partial x^2}\right|_{i,j}^{L} = a_i \varphi_{i-1,j} + b_i \varphi_{i,j} + c_i \varphi_{i+1,j} \tag{2}$$

where:

$$a_i = \frac{2}{\Delta x_i(\Delta x_{i+1}+\Delta x_i)}, \quad b_i = -\frac{2}{\Delta x_i \Delta x_{i+1}}, \quad c_i = \frac{2}{\Delta x_{i+1}(\Delta x_{i+1}+\Delta x_i)} \tag{3}$$

For a non-uniform grid, the higher order compact approximation for the second derivative in the $x$-direction may be rewritten in the following general form:

$$\alpha_i \frac{\partial^2 \varphi}{\partial x^2}\bigg|_{i-1,j}^{H} + \frac{\partial^2 \varphi}{\partial x^2}\bigg|_{i,j}^{H} + \beta_i \frac{\partial^2 \varphi}{\partial x^2}\bigg|_{i+1,j}^{H} = a_i \varphi_{i-1,j} + b_i \varphi_{i,j} + c_i \varphi_{i+1,j} \tag{4}$$

where the coefficients $\alpha_i$, $\beta_i$, $a$, $b$, and $c$ are functions of the non-uniform grid spacing $\Delta x_i$. Matching the Taylor series of various orders can derive the relations between the coefficients in Eq. (4). The coefficients used in Eq. (4) are:

$$\alpha_i = \frac{\Delta x_{i+1}(\Delta x_i^2 + \Delta x_i \Delta x_{i+1} - \Delta x_{i+1}^2)}{\Delta x_i^3 + \Delta x_{i+1}^3 + 4\Delta x_{i+1}\Delta x_i^2 + 4\Delta x_i \Delta x_{i+1}^2}, \quad \beta_i = \frac{\Delta x_i(\Delta x_{i+1}^2 + \Delta x_i \Delta x_{i+1} - \Delta x_i^2)}{\Delta x_i^3 + \Delta x_{i+1}^3 + 4\Delta x_{i+1}\Delta x_i^2 + 4\Delta x_i \Delta x_{i+1}^2},$$

$$a_i = \frac{12\Delta x_{i+1}}{\Delta x_i^3 + \Delta x_{i+1}^3 + 4\Delta x_{i+1}\Delta x_i^2 + 4\Delta x_i \Delta x_{i+1}^2}, \quad b_i = \frac{12}{\Delta x_{i+1}^2 + 3\Delta x_{i+1}\Delta x_i + \Delta x_i^2},$$

$$c_i = \frac{12\Delta x_i}{\Delta x_i^3 + \Delta x_{i+1}^3 + 4\Delta x_{i+1}\Delta x_i^2 + 4\Delta x_i \Delta x_{i+1}^2} \tag{5}$$

For non-periodic boundaries, Eq. (4) can no longer be applied to the boundary points so that boundary schemes are required. The second order derivative approximation at the boundary point in the x-direction can be written as:

$$\frac{\partial^2 \varphi}{\partial x^2}\bigg|_{i,j}^{H} + \beta_i \frac{\partial^2 \varphi}{\partial x^2}\bigg|_{i+1,j}^{H} = a_i \varphi_{i,j} + b_i \varphi_{i+1,j} + c_i \varphi_{i+2,j} \tag{6}$$

The coefficients in Eq. (6) can be calculated by matching the Taylor series of various orders as follows:

$$\beta_i = \frac{2\Delta x_{i+1} + \Delta x_{i+2}}{\Delta x_{i+1} - \Delta x_{i+2}}, \quad a_i = \frac{6}{\Delta x_{i+1}^2 - \Delta x_{i+2}^2},$$

$$b_i = \frac{6}{(\Delta x_{i+2} - \Delta x_{i+1})\Delta x_{i+2}}, \quad c_i = \frac{6}{\Delta x_{i+1}\Delta x_{i+2}(\Delta x_{i+1}^2 - \Delta x_{i+2}^2)} \tag{7}$$

The same approximations with non-uniform space $\Delta y_j$ are also applied to the spatial derivative terms in the y-direction. When the value of $\varphi_{i,j}$ is given, the value of second order derivative at every point can be calculated by solving Eq. (4) and Eq. (6). The tridiagonal matrix algorithm (TDMA) are applied for solving the equations.

**2.2 Compact correction term on non-uniform grid**

For classical finite difference scheme, the derivative for every point can be expressed by the linear combination of the point function in Eq. (2), which means that the classical finite difference equation for Poisson's equation can be directly derived from the difference approximations of the derivative. The classical finite difference formulation with lower accuracy can be written as follows:

$$\frac{\partial^2 \varphi}{\partial x^2}\bigg|_{i,j}^{L} + \frac{\partial^2 \varphi}{\partial y^2}\bigg|_{i,j}^{L} = f(x_i, y_j) \tag{8}$$

The resulting linearized equations can be derived by substituting Eq. (2) into Eq. (8) as follows:

$$a_p \varphi_{i,j} = a_w \varphi_{i-1,j} + a_e \varphi_{i+1,j} + a_s \varphi_{i,j-1} + a_n \varphi_{i,j+1} + b \tag{9}$$

The coefficients are functions of the non-uniform grid spacing $\Delta x_i$ and $\Delta y_j$, which can be written as:

$$a_w = \frac{2}{\Delta x_i(\Delta x_i + \Delta x_{i+1})}, \quad a_e = \frac{2}{\Delta x_{i+1}(\Delta x_i + \Delta x_{i+1})}, \quad a_s = \frac{2}{\Delta y_j(\Delta y_j + \Delta y_{j+1})},$$

$$a_n = \frac{2}{\Delta y_{j+1}(\Delta y_j + \Delta y_{j+1})}, \quad a_p = \frac{2}{\Delta x_i \Delta x_{i+1}} + \frac{2}{\Delta y_j \Delta y_{j+1}}, \quad b = -f(x_i, y_j) \tag{10}$$

The numerical solution $\varphi_{i,j}^*$ with lower accuracy of Poisson's equation can be obtained by solving Eq. (9) using the Gauss-Seidel method. The derivations of classical finite difference formulation are simple and straightforward because of its explicit form. Compared with the classical schemes, compact schemes with the same stencil width are implicit schemes with more accuracy in Eq. (4) and in Eq. (6). The derivations of compact finite difference formulation are not straightforward. The existence of implicit difference approximation creates extra complexity in deriving the compact finite difference equation. The processes of matrix transformation are inflexible and time consuming. Therefore, a new idea by combining compact and classical finite difference formulations are proposed for solving Poisson's equation, which can improve the flexibility of derivation process and increase the accuracy of numerical solution. The modified finite difference formulation from Eq. (9) can be written as:

$$a_p \varphi_{i,j} = a_w \varphi_{i-1,j} + a_e \varphi_{i+1,j} + a_s \varphi_{i,j-1} + a_n \varphi_{i,j+1} + b + b^* \tag{11}$$

The coefficients of the above equation can be written as:

$$a_w = \frac{2}{\Delta x_i(\Delta x_i + \Delta x_{i+1})}, \quad a_e = \frac{2}{\Delta x_{i+1}(\Delta x_i + \Delta x_{i+1})}, \quad a_s = \frac{2}{\Delta y_j(\Delta y_j + \Delta y_{j+1})},$$

$$a_n = \frac{2}{\Delta y_{j+1}(\Delta y_j + \Delta y_{j+1})}, \quad a_p = \frac{2}{\Delta x_i \Delta x_{i+1}} + \frac{2}{\Delta y_j \Delta y_{j+1}}, \quad b = -f(x_i, y_j) \tag{12}$$

where the source term $b^*$ is the compact correction term, which is designed for increasing the accuracy of numerical solution. The compact correction term based on compact and classical finite difference formulation for second derivative can be written as:

$$b^* = \left.\frac{\partial^2 \varphi}{\partial x^2}\right|_{i,j}^H + \left.\frac{\partial^2 \varphi}{\partial y^2}\right|_{i,j}^H - \left.\frac{\partial^2 \varphi}{\partial x^2}\right|_{i,j}^L + \left.\frac{\partial^2 \varphi}{\partial y^2}\right|_{i,j}^L \tag{13}$$

When the value of $\varphi_{i,j}$ is given, the compact correction term can be calculated by solving the classical formulation Eq. (2) and compact difference formulation Eqs. (4) and (6). In fact, the exact solution of $\varphi_{i,j}$ can be not given to calculate the value of compact correction term. However, the value of this term can be estimated by numerical solutions $\varphi_{i,j}^*$ from Eq. (9). Based on the estimated value of compact correction term, the modified finite difference formulation Eq. (11) can be applied to solve Poisson's equation. The improved finite difference method based on compact correction term can be developed using this approach.

### 2.3 Poisson's equation solving process
The Poisson's equations mentioned above may be solved numerically and solution was marched according to the following steps:

Step 1: Calculate the temporary numerical solutions $\varphi_{i,j}^*$ from Eq. (9)

Step 2: According to the temporary numerical solution $\varphi_{i,j}^*$, calculate the values of second derivative $\left.\frac{\partial^2 \varphi}{\partial x^2}\right|_{i,j}^L$, $\left.\frac{\partial^2 \varphi}{\partial y^2}\right|_{i,j}^L$ by solving Eq. (2) and $\left.\frac{\partial^2 \varphi}{\partial x^2}\right|_{i,j}^H$, $\left.\frac{\partial^2 \varphi}{\partial y^2}\right|_{i,j}^H$ by solving Eq. (4), Eq. (6), respectively.

Step 3: According to the value of second derivative $\left.\dfrac{\partial^2 \varphi}{\partial x^2}\right|_{i,j}^{L}$, $\left.\dfrac{\partial^2 \varphi}{\partial y^2}\right|_{i,j}^{L}$, $\left.\dfrac{\partial^2 \varphi}{\partial x^2}\right|_{i,j}^{H}$ and $\left.\dfrac{\partial^2 \varphi}{\partial y^2}\right|_{i,j}^{H}$, calculate the value of compact correction term $b^*$ by Eq. (13).

Step 4: Calculate the numerical solution $\varphi_i$ from Eq. (11)

## 2.4 Extension to three-dimensional Poisson's equation

The present method has higher accuracy and also provides a direct extension to the three-dimensional Poisson's equation, even to the multi-dimensional cases. The 3D Poisson's equation can be written as

$$\frac{\partial^2 u}{\partial x^2}+\frac{\partial^2 u}{\partial y^2}+\frac{\partial^2 u}{\partial z^2}=f(x,y,z) \quad (x,y,z)\in \Omega \tag{14}$$

where $\Omega = [0, L_x] \times [0, L_y] \times [0, L_z]$ is a cubic domain with Dirichlet boundary condition. The spatial domain $\Omega$ can be discretized with non-uniform grid sizes $\Delta x_i = x_{i+1} - x_i$, $\Delta y_j = y_{j+1} - y_j$ and $\Delta z_k = z_{k+1} - z_k$, respectively. The three-dimensional finite difference formulation based on compact correction term can be written as:

$$a_p \varphi_{i,j,k} = a_w \varphi_{i-1,j,k} + a_e \varphi_{i+1,j,k} + a_s \varphi_{i,j-1,k} + a_n \varphi_{i,j+1,k} + a_f \varphi_{i,j,k-1} + a_b \varphi_{i,j,k+1} + b + b^* \tag{15}$$

The coefficients of Eq. (15) can be written as:

$$a_w = \frac{2}{\Delta x_i (\Delta x_i + \Delta x_{i+1})}, \quad a_e = \frac{2}{\Delta x_{i+1}(\Delta x_i + \Delta x_{i+1})}, \quad a_s = \frac{2}{\Delta y_j (\Delta y_j + \Delta y_{j+1})},$$

$$a_n = \frac{2}{\Delta y_{j+1}(\Delta y_j + \Delta y_{j+1})}, \quad a_f = \frac{2}{\Delta z_k (\Delta z_k + \Delta z_{k+1})}, \quad a_b = \frac{2}{\Delta z_{k+1}(\Delta z_k + \Delta z_{k+1})},$$

$$a_p = \frac{2}{\Delta x_i \Delta x_{i+1}} + \frac{2}{\Delta y_j \Delta y_{j+1}} + \frac{2}{\Delta z_k \Delta z_{k+1}}, \quad b = -f(x_i, y_j, z_k),$$

$$b^* = \left.\frac{\partial^2 \varphi}{\partial x^2}\right|_{i,j,k}^{H} + \left.\frac{\partial^2 \varphi}{\partial y^2}\right|_{i,j,k}^{H} + \left.\frac{\partial^2 \varphi}{\partial z^2}\right|_{i,j,k}^{H} - \left.\frac{\partial^2 \varphi}{\partial x^2}\right|_{i,j,k}^{L} - \left.\frac{\partial^2 \varphi}{\partial y^2}\right|_{i,j,k}^{L} - \left.\frac{\partial^2 \varphi}{\partial z^2}\right|_{i,j,k}^{L} \tag{16}$$

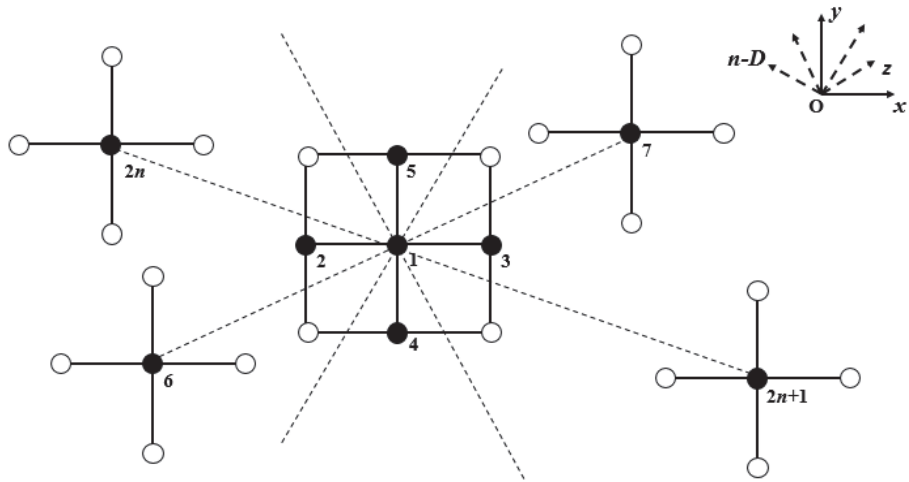

Figure 1 2n+1 points on a stencil of CCFDM

When the compact correction term $b^*$ is equal to zero, Eq. (15) represents the discretization formulation of classical finite difference method with lower accuracy. Firstly, the numerical solution $\varphi^*_{i,j,k}$ can be obtained by classical finite difference method. Secondly, the estimated value of $b^*$ can be calculated by solving Eq. (16). Finally, the numerical solution with higher order accuracy can be obtained by solving Eq. (15). Based on the proposed discretization methodology, the improved finite difference method with compact correction term can be easily developed to the multi-dimensional cases by the similar deduce process. Therefore, the CCFDM can be applied not only for solving two-dimensional Poisson's equation, but can also be applied to solve three-dimensional Poisson's equation on non-uniform grid. For two- and three-dimensional case, the discretized equations are established on a stencil of five points and seven points, respectively. The $n$-dimensional discretization of Poisson's equation can be developed on a stencil of 2n+1 points, as is shown in Fig. 1.

## 3 NUMERICAL EXPERIMENTS

In this section, the numerical results for four examples with smooth and finite regular solutions will be presented using the proposed method. The Poisson's equation can be solved with Dirichlet boundary condition. Our code is written in Fortran 77 on a personal laptop equipped with Intel(R) Core (TM) 2 Duo CPU (2.40Hz) and 2.00 GB RAM. The initial value of solutions are set to the zero vector. The iteration stops when the average value of the residual vector is reduced by $10^{-14}$. The maximum absolute errors $e_{max}$ over the discretized grid are considered, and the order of convergence of the method is evaluated by the following formula.

$$\text{order} = \frac{\log(e_1 / e_2)}{\log 2} \quad (17)$$

where $e_1$ and $e_2$ are the errors of numerical solutions for two grid systems with different numbers of grid points. The grid spacing for the first grid system is twice as larger than that for the second gird system. In this paper, the maximum absolute errors $e_{max}$ and the mean absolute errors $e_{ave}$ are applied for analyzing the characteristics of both FDM and CCFDM. The value of $e_{ave}$ is equal to the averaged value of the errors between the exact and numerical solutions. To investigate the influence of grid system on the numerical solutions, three different types of grid systems are used for solving Poisson's equations. They have different stretch ratio of the neighbor grid, which include uniform grid, sinh-based grid defined by

$$x_i = L_x \frac{\sinh(\gamma \xi_i)}{\sinh(\gamma)} \quad \text{and} \quad y_j = L_y \frac{\sinh(\gamma \xi_j)}{\sinh(\gamma)} \quad (18)$$

and tanh-based grid defined by

$$x_i = L_x \frac{\tanh(\gamma \xi_i)}{\tanh(\gamma)} \quad \text{and} \quad y_j = L_y \frac{\tanh(\gamma \xi_j)}{\tanh(\gamma)} \quad (19)$$

where γ is the control parameter of grid systems. The grid intervals can be revised by changing the value of γ. $\xi$ is the fundamental grid points with uniform intervals in the domain [-1, 1]. $L_x$ and $L_y$ are half the interval length in the $x$- and $y$- directions, respectively. For the uniform mesh, the grid interval is equal to $\Delta x = L_x/M$ or $\Delta y = L_y/N$ in the $x$- and $y$- coordinate directions, respectively.

The first type of the grid is one of the most commonly used grid system, its stretch ratio of the neighbor grid is 1.0. For the second type of grid, the density of grid points become more and higher along the positive direction of axes. The third one is contrary to the second one.

### 3.1 Problem 1
Consider the following Poisson's equation:

$$\frac{\partial^2 \varphi}{\partial x^2}+\frac{\partial^2 \varphi}{\partial y^2}=-2\pi^2\sin(\pi x)\cos(\pi y) \quad (x,y)\in \Omega=[0,1]\times[0,1] \tag{20}$$

The exact solution of Eq. (20) is
$$\varphi(x,y)=\sin(\pi x)\cos(\pi y) \tag{21}$$

The above two-dimensional Poisson's equation was solved on uniform grid to validate the improved finite difference method with compact correction term. Table 1 shows the maximum absolute errors, the average absolute errors and other information of classical finite difference method and compact correction method. The numerical solutions obtained by the improved finite difference method have higher accuracy than that of classical finite difference scheme obviously. The results indicate that the FDM based on a stencil of five points can obtain the numerical solutions with second order accuracy, while the numerical solutions obtained by the CCFDM has four order accuracy based on the same stencil. The improved finite difference method need higher computational cost (CPU time) than that of the classical finite difference method with the same number of grid. That is because the compact correction term needs additional computational time in the calculating process of CCFDM. When the number of grid is 10×10, the CPU time can be neglected because the CPU time is less than 0.02. The CPU cost tends to dramatically escalate with the increasing number of grid points. The data in Table 1 shows that the cost is over 55 CPU seconds with the number of grid 160×160. The values of maximum and average absolute errors from the FDM are equal to $1.08\times10^{-5}$ and $4.53\times10^{-5}$, respectively. For the improved finite difference method, more accurate numerical solution can be obtained with the number of grid 40×40. The corresponding data in Table 1 shows that the computational cost of CCFDM is less than one CPU second. Therefore, the improved finite difference method with the fourth order accuracy is more than 50 times faster than the classical finite difference method with the second order accuracy.

Table 1 Errors, average absolute error, order of convergence and CPU seconds of both FDM and CCFDM for solving the Problem 1.

| grid | FDM | | | | | CCFDM | | | | |
|---|---|---|---|---|---|---|---|---|---|---|
| | $e_{max}$ | order | $e_{ave}$ | order | CPU | $e_{max}$ | order | $e_{ave}$ | order | CPU |
| 10×10 | 2.76(-3) | --- | 1.35(-3) | --- | --- | 5.19(-4) | --- | 2.01(-4) | --- | --- |
| 20×20 | 6.91(-4) | 2.0 | 3.14(-4) | 2.1 | 0.02 | 3.91(-5) | 3.7 | 1.22(-5) | 4.0 | 0.05 |
| 40×40 | 1.73(-4) | 2.0 | 7.52(-5) | 2.1 | 0.23 | 2.66(-6) | 3.9 | 7.32(-7) | 4.1 | 0.40 |
| 80×80 | 4.33(-5) | 2.0 | 1.84(-5) | 2.0 | 3.18 | 1.72(-7) | 4.0 | 4.08(-8) | 4.2 | 5.13 |
| 160×160 | 1.08(-5) | 2.0 | 4.53(-6) | 2.0 | 55.48 | 1.06(-8) | 4.0 | 2.54(-9) | 4.0 | 85.06 |

**3.2 Problem 2**
Consider the following Poisson's equation
$$\frac{\partial^2 \varphi}{\partial x^2}+\frac{\partial^2 \varphi}{\partial y^2}=2e^{x+y} \quad (x,y)\in \Omega=[0,1]\times[0,1] \tag{22}$$

The exact solution of Eq. (22) is
$$\varphi(x,y)=e^{x+y} \tag{23}$$

Table 2 shows the maximum absolute errors and the convergence order on three different kinds of grid systems. The value of the control parameter $\gamma$ used here is equal to 1.0 for both sinh-based grid and tanh-based grid. With the same number of grid points, the different stretch ratio of the neighbor grid can influence the accuracy of numerical solutions. For both FDM and CCFDM, the numerical solutions obtained based on Sinh-based grid system have higher accuracy than the solutions based on uniform and Tanh-based grid systems. Because the gradient value of exact solutions become larger gradually along the positive direction of axes. For the Sinh-based grid system, the grid density become higher along the same direction. The ratio of the neighbor grid for Sinh-based grid system can adapt to the change of exact solutions better. Correspondingly, the change trend of grid density for Tanh-based grid system are not consistent with the exact solutions completely. The numerical solutions based on Tanh-based grid system have lower accuracy. The data in Table 2 indicates that the convergence order of classical finite difference method are almost fixed to the second order. However, the convergence order of CCFDM based on those data oscillated nearby the fourth order slightly. For most cases, the convergence order of numerical solutions obtained by the CCFDM is lower than the fourth order accuracy, especially for the case of fewer grid points. The reason for this phenomenon is related to the calculation of compact correction term. The compact correction term can be obtained based on the estimated solutions, instead of exact solutions. The influence of the estimated solutions on the numerical results can be seen from the convergence order when the number of grid points is very small. However, the influence of estimated solutions on the final results are very limited, and it can be neglected due to the dramatic increase of accuracy. The data in this table indicate that the convergence order of improved finite difference method with compact correction term is close to the fourth order.

Figure 2 shows the change of the average absolute errors on Sinh-based grid system with the different values of control parameters $\gamma$. When the value of $\gamma$ approach zero infinitely, the Sinh-based grid system with the uniform grid density can be considered as uniform grid system. With the increasing value of $\gamma$, the stretch ratio of the neighbor grid increases and the change of grid density become more obvious. The average absolute errors of both FDM and CCFDM decreases at first and then increase with the control parameters. The minimum value of the average absolute errors can be achieved at about $\gamma=0.75$ for the FDM and $\gamma=0.6$ for the CCFDM. When the stretch ratio of neighbor grid changes, the fluctuations of the average errors for the FDM are more apparent than that for the CCFDM. In the range of control parameter $\gamma$ from 0.01 to 1, the maximum difference of the average absolute errors is equal to $9.63\times10^{-6}$ for the FDM and $6.57\times10^{-9}$ for the CCFDM, respectively. Therefore, the change of grid system has more influence on the numerical solutions of the FDM than that of the CCFDM. The improved finite difference method with compact correction term can be beneficial to reduce the dependence of numerical solutions on grid systems.

Table 2 The maximum absolute errors and the convergence order of FDM and CCFDM on different grid systems for solving Problem 2.

|  | FDM | | | CCFDM | | |
| --- | --- | --- | --- | --- | --- | --- |
|  | Sinh-based grid | Uniform grid | Tanh-based grid | Sinh-based grid | Uniform grid | Tanh-based grid |
| 10×10 | 3.046(-4) | 3.549(-4) | 1.465(-3) | 3.535(-5) | 4.486(-5) | 1.210(-4) |
| order | --- | --- | --- | --- | --- | --- |
| 20×20 | 7.659(-5) | 8.941(-5) | 3.748(-4) | 2.587(-6) | 3.198(-6) | 9.085(-6) |
| order | 2.0 | 2.0 | 2.0 | 3.8 | 3.8 | 3.9 |
| 40×40 | 1.917(-5) | 2.248(-5) | 9.391(-5) | 1.815(-7) | 2.140(-7) | 6.220(-7) |
| order | 2.0 | 2.0 | 2.0 | 3.8 | 3.9 | 3.9 |
| 80×80 | 4.795(-6) | 5.624(-6) | 2.351(-5) | 1.232(-8) | 1.385(-8) | 4.070(-8) |
| order | 2.0 | 2.0 | 2.0 | 3.9 | 3.9 | 3.9 |
| 160×160 | 1.199(-6) | 1.406(-6) | 5.879(-6) | 8.156(-10) | 8.808(-10) | 2.603(−9) |
| order | 2.0 | 2.0 | 2.0 | 3.9 | 4.0 | 4.0 |

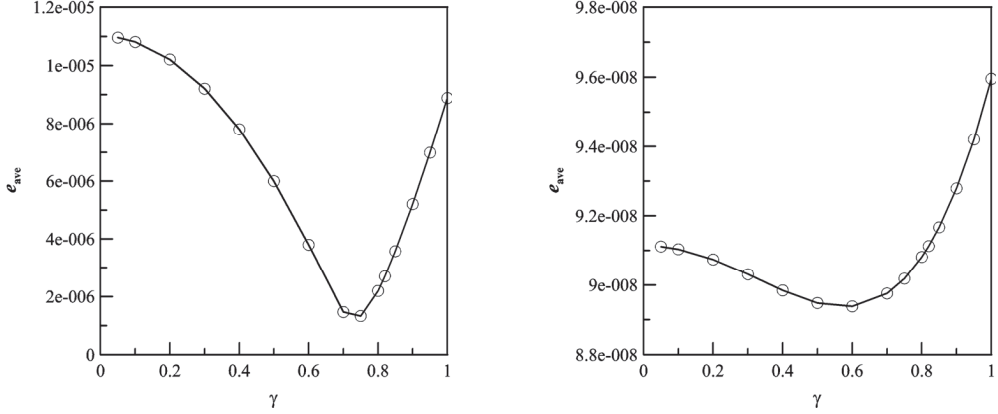
Figure 2 Comparison of the average absolute errors of FDM and CCFDM on Sinh-based grid system with different control parameters.

### 3.3 Problem 3
Consider the following Poisson's equation
$$\frac{\partial^2 \varphi}{\partial x^2}+\frac{\partial^2 \varphi}{\partial y^2}+\frac{\partial^2 \varphi}{\partial z^2}=(8\pi^2-1)\times e^{-2\pi x-2\pi y}\sin(z) \quad (x,y,z)\in \Omega=[0,1]\times[0,1]\times[0,1] \qquad (24)$$
The analytic solution of Eq. (24) is
$$\varphi(x,y,z)=e^{(-2\pi x-2\pi y)}\sin(z) \qquad (25)$$

The errors and convergence order of CCFDM for solving three dimensional Poisson's equations on uniform and non-uniform grid are summarized in Table 3. The numerical results showed the advantage of using the suitable non-uniform grid system for the CCFDM. The non-uniform grid system used here is the Tanh-based grid system with the control parameter of $\gamma=1.1$ along the $x$- and $y$- directions. For the case of fewer grid points, the CPU time for using non-uniform grid are more than that for using uniform grid because of the grid construction. The accuracy of numerical solutions increases with the increasing numbers of grid points. Meanwhile, it leads to the obvious increase of the CPU time and iterations. When the number of grid is over than 20×20×20, the CCFDM of using uniform grid takes more CPU seconds due to more number of iterations. When the number of grid points used in calculation is small enough, it is inevitable for the CCFDM to deteriorate the convergence order. It can be observed in table 3 that the convergence rate on uniform grid decreases to 2.8th order when the number of grid decreases to 20×20×20, while the convergence rate on non-uniform grid can still achieve 3.5th order at the same parameters.

Table 3 Errors, convergence order and iterations on uniform and non-uniform grid for solving Problem 3.

| Grid | $e_{max}$ | order | $e_{ave}$ | order | CPU | iteration |
|---|---|---|---|---|---|---|
| Uniform grid system | | | | | | |
| 10×10×10 | 1.80(-3) | --- | 1.16(-4) | --- | 0.02 | 440 |
| 20×20×20 | 2.53(-4) | 2.8 | 1.01(-5) | 3.5 | 1.07 | 1638 |
| 40×40×40 | 2.46(-5) | 3.4 | 7.43(-7) | 3.8 | 45.08 | 5954 |
| 80×80×80 | 1.95(-6) | 3.7 | 5.02(-8) | 3.9 | 2777.94 | 21251 |
| Non-uniform grid system | | | | | | |
| 10×10×10 | 3.02(-4) | --- | 2.97(-5) | --- | 0.03 | 393 |
| 20×20×20 | 2.63(-5) | 3.5 | 1.82(-6) | 4.0 | 0.98 | 1458 |
| 40×40×40 | 1.94(-6) | 3.8 | 1.11(-7) | 4.0 | 39.66 | 5260 |
| 80×80×80 | 1.32(-7) | 3.9 | 6.83(-9) | 4.0 | 2544.86 | 18627 |

The maximum error occurs near the domain of computational bound where the gradient of exact solutions change obviously. The non-uniform grid system can produce higher grid density at this domain than uniform grid system with the same number of grid points. By this way, the CCFDM using non-uniform grid can reduce the maximum absolute errors efficiently. The maximum absolute error is $1.95\times10^{-6}$ with $80\times80\times80$ grid on uniform grid, while almost the same accuracy can be achieved with $40\times40\times40$ grid on non-uniform grid. Therefore, the numerical solutions on non-uniform grid can achieve significantly better accuracy than those on uniform grid and be obtained with the shorter computational time.

### 3.4 Problem 4
Consider the following Poisson's equation
$$\frac{\partial^2 \varphi}{\partial x^2}+\frac{\partial^2 \varphi}{\partial y^2}+\frac{\partial^2 \varphi}{\partial z^2}+\frac{\partial^2 \varphi}{\partial u^2}=4e^{x+y+z+u} \quad (x,y,z,u)\in \Omega=[0,1]\times[0,1]\times[0,1]\times[0,1]\times[0,1] \quad (26)$$
The exact solution of Eq. (26) is
$$\varphi(x,y)=e^{x+y+z+u+v} \quad (27)$$

Table 4 Errors, convergence order, CPU time and iterations of the FDM and CCFDM for solving Problem 4.

| FDM | | | | | | |
|---|---|---|---|---|---|---|
| Grid | $e_{max}$ | order | $e_{ave}$ | order | CPU | iteration |
| 10×10×10×10 | 1.07(-3) | | 4.16(-4) | | 0.37 | 321 |
| 20×20×20×20 | 2.77(-4) | 1.9 | 8.95(-5) | 2.2 | 42.83 | 1245 |
| 30×30×30×30 | 1.23(-4) | 2.0 | 3.76(-5) | 2.2 | 693.58 | 2733 |
| 40×40×40×40 | 6.94(-5) | 2.0 | 2.05(-5) | 2.1 | 3978.13 | 4768 |
| 60×60×60×60 | 3.09(-5) | 2.0 | 8.84(-6) | 2.1 | 46878.35 | 10434 |
| CCFDM | | | | | | |
| Grid | $e_{max}$ | order | $e_{ave}$ | order | CPU | |
| 10×10×10×10 | 1.55(-4) | | 6.78(-5) | | 0.53 | 542 |
| 20×20×20×20 | 1.19(-5) | 3.7 | 4.23(-6) | 4.0 | 56.37 | 2025 |
| 30×30×30×30 | 2.70(-6) | 3.7 | 8.35(-7) | 4.0 | 937.89 | 4344 |
| 40×40×40×40 | 9.50(-7) | 3.6 | 2.64(-7) | 4.0 | 5288.13 | 7446 |
| 60×60×60×60 | 2.12(-7) | 3.7 | 5.23(-8) | 4.0 | 62313.80 | 15861 |

The improved finite difference method with compact correction term can be used for multi-dimensional Poisson's equations. Four-dimensional Poisson's equation as a typical example of multidimensional problem are solved to verify the proposed method. The non-uniform grid system for solving Problem 4 is the sinh-based grid system with the control parameters $\gamma=1$. Table 4 shows that the errors, CPU time and iterations of FDM and CCFDM on non-uniform grid. Compared with the calculations of low dimensional problems, the CPU time and the number of iterations for solving multi-dimensional problem increase obviously. That is because the significant increase of grid points from two-dimensional case to multi-dimensional case. The accuracy of numerical solutions can be improved by the finer grid system while the CPU time and the number of iterations increase dramatically. The CPU seconds are more than doubled when the number of grid points are doubled. It is necessary for using an effective method to obtain higher order numerical solutions with the coarse grid. The CCFDM need to take more CPU seconds and iterations than the FDM with the same parameters. However, the accuracy of numerical solutions for the CCFDM is higher than that for the FDM. The numerical results demonstrate that the classical finite difference method generates the second order accuracy while the convergence rate of the proposed CCFDM approach to four order. The advantage of the CCFDM become more obvious for solving multi-dimensional Poisson's equations. The maximum and average absolute errors of the FDM with 60×60×60×60 grid are 3.09×10$^{-5}$ and 8.84×10$^{-6}$, respectively. The data in Table 4 show the cost of the FDM is over 46,878 CPU seconds with 10,434 steps of iterations. The numerical solutions with higher accuracy can be obtained by the CCFDM using 20×20×20×20 grid. The computational cost of the CCFDM is only 42.83 CPU seconds, which is about 1,000 times faster than the FDM. The similar conclusions can be obtained by the comparisons using other data. The classical finite difference method with second order accuracy is not competitive due to its lower accuracy, the extra computations for the improved finite difference method with compact correction term are negligible compared to the increased accuracy. Therefore, the improved finite difference method with compact correction term provide an efficient way to solve the Poisson's equations.

## 4 CONCLUSIONS

An improved finite difference method with compact correction term (CCFDM) on non-uniform grid are developed for solving the two to multi-dimensional Poisson's equations. The compact correction term is established by the coupled high order compact and low order classical finite difference formulations. The proposed method is developed on a stencil of 2n+1 points and it can approach to the four order accuracy for solving N-dimensional Poisson's equation.

The results indicate that the classical finite difference method with the second order accuracy is not competitive due to its low accuracy, the extra computations for the proposed method are negligible compared to the increased accuracy. When the stretch ratio of neighbor grid changes, the CCFDM is beneficial to reduce the dependence of numerical solutions on grid systems because of the leading role of the implicit compact finite difference formulation in calculation. The numerical solutions on non-uniform grid can achieve significantly better accuracy than those on uniform grid while the CPU seconds and the number of iterations using non-uniform grid are less than that using uniform grid with the same parameters. The computational time of both FDM and CCFDM for solving from two to multi-dimensional Poisson's equations increases dramatically. The advantage of the CCFDM with higher order accuracy become obvious for solving multi-dimensional problems. The basic idea of adding compact correction term into finite difference method in this paper can also be extended to solve other partial difference equations, such as 2D and 3D convection-diffusion equations.

## ACKNOWLEDGEMENT

The financial supports from the Chinese National Natural Science Foundation under Grants No. 50876067 and China Postdoctoral Science Foundation funded project No. 2014M562475 are gratefully acknowledged.

## NOMENCLATURE

$a_e$    coefficients of discretized Poisson's equation

$a_i$    coefficients of discretized equation for second derivative

$a_n$    coefficients of discretized Poisson's equation

$a_p$    coefficients of discretized Poisson's equation

$a_s$    coefficients of discretized Poisson's equation

$b_i$    coefficients of discretized equation for second derivative

$b^*$    compact correction term

$c_i$    coefficients of discretized equation for second derivative

$f$    forcing function

$H$    numerical method with higher order accuracy

*L*    numerical method with lower order accuracy
$L_x$    horizontal length of computational domain
$L_y$    vertical length of computational domain
*M*    number of grid points in *x*-direction
*N*    number of grid points in *y*-direction
*x*    horizontal coordinate
$x_i$    grid point in x-direction
*y*    vertical coordinate
$y_i$    grid point in y-direction

**Greek Symbols**

$\alpha_i$    coefficient of discretized equations for second derivative
$\beta_i$    coefficient of discretized equations for second derivative
$\gamma$    stretch ratio of neighbor grid
$\xi$    fundamental grid points with uniform intervals in the domain [-1, 1]
$\Delta x_i$    grid spacing in *x*-direction
$\Delta y_i$    grid spacing in *y*-direction
$\varphi$    variable value
$\varphi_{i,j}$    numerical solution
$\varphi_{i,j}^*$    numerical solution with lower accuracy
$\Omega$    computational domain